\documentclass{amsart}

\usepackage{amsmath,amsthm, amssymb}

\usepackage[all]{xy}
\usepackage{graphicx}
\usepackage{amsmath}
\usepackage{amsfonts}

\newtheorem{theorem}{Theorem}[section]

\theoremstyle{remark}

\newtheorem{rmk}[theorem]{Remark}
\newtheorem{ex}[theorem]{Example}

\numberwithin{equation}{section}

\newcommand{\lap}{\vartriangle}

\newcommand{\R}{\mathbb{R}}

\newcommand{\pr}{\Pi_\lap}

\newcommand{\lef}{L\!e\!f}
\newcommand{\rel}{\mathcal{R}}

\newcommand{\eps}{\epsilon}
\newcommand{\epseta}{\eps_\eta}
\newcommand{\ixi}{i_\xi}
\DeclareMathOperator{\Ext}{Ext}

\newcommand{\Rr}{\mathbb R}
\newcommand{\Z}{\mathbb{Z}}
\newcommand{\T}{\mathbb{T}}

\renewcommand{\d}{\mathrm d}
\newcommand{\g}{\mathfrak{g}}

\begin{document}

\title{Examples of compact $K$-contact manifolds with no Sasakian metric}

\author[B. Cappelletti-Montano]{Beniamino Cappelletti-Montano}
 \address{Dipartimento di Matematica e Informatica, Universit\`{a} degli Studi
 di Cagliari, Via Ospedale 72, 09124 Cagliari, Italy}
 \email{b.cappellettimontano@gmail.com}

 \author[A. De Nicola]{Antonio De Nicola}
 \address{CMUC, Department of Mathematics, University of Coimbra, 3001-501 Coimbra, Portugal}
 \email{antondenicola@gmail.com}

 \author[J.~C. Marrero]{Juan Carlos Marrero}
 \address{Unidad Asociada ULL-CSIC ``Geometr{\'\i}a Diferencial y Mec\'anica Geo\-m\'e\-tri\-ca''
Departamento de Matem\'aticas, Estad{\'\i}stica e Investigaci\'on Operativa, Facultad de Ciencias, Universidad de La Laguna, La Laguna, Tenerife, Spain}
 \email{jcmarrer@ull.es}

\author[I. Yudin]{Ivan Yudin}
 \address{CMUC, Department of Mathematics, University of Coimbra, 3001-501 Coimbra, Portugal}
 \email{yudin@mat.uc.pt}

\maketitle

\begin{abstract}
Using the Hard Lefschetz Theorem for Sasakian manifolds, we find two examples of compact $K$-contact nilmanifolds with no compatible Sasakian metric in dimensions five and seven, respectively.
\end{abstract}


\section{Introduction}
Construction of examples of compact symplectic manifolds with no K\"{a}hler structures is a topic which attracted wide interest in recent years (see e.g. \cite{tralleoprea1997} and reference therein).  Many techniques have been used to identify such examples. The first technique used to prove that some examples were non-K\"{a}hler consists in showing that some odd Betti number of the manifold is not even and thus the manifold cannot be K\"{a}hler (see e.g. \cite{thurston}). In other cases (e.g. \cite{bensongordon}), when the all odd Betti numbers are even, one possibility to prove that a certain symplectic manifold cannot admit a K\"{a}hler structure is to show that it does not satisfy the Hard Lefschetz theorem.

In odd dimensions, a Hard Lefschetz theorem for compact co-K\"{a}hler manifolds was proven in \cite{chineadeleonmarrero}. However, until recently only the Betti number technique was available to show that a compact $K$-contact manifold does not admit any Sasakian structure.
An example of a compact $K$-contact manifold with no Sasakian structure was found by Boyer and Galicki in \cite[Example 7.4.16]{boyergalicki2008}.
They constructed a non-trivial $\T^3$-bundle over $\T^2$ and proved that it is non-Sasakian since the first Betti number is odd.
Recently a Hard Lefschetz theorem for compact Sasakian manifolds was proven (\cite{cappellettidenicolayudin2013}).

In this paper we give two examples of compact co-orientable $K$-contact manifolds with no Sasakian structure compatible with the contact one. The examples, in dimensions five and seven respectively, are intended as first applications of our Hard Lefschetz theorem. We take advantage of the fact that lately Kutsak has classified the invariant contact structures  on nilmanifolds up to dimension seven in \cite{kutsak2012}. Thus, by using this classification one can construct examples of compact contact nilmanifolds in any odd dimension $3$, $5$, $7$. Actually, since the only $3$-dimensional nilpotent Lie algebra admitting an invariant contact structure is the Heisenberg algebra, the  non-trivial dimensions are $5$ and $7$. Among these examples we find one in dimension $5$ with $b_1 = 2$ and one in dimension $7$ with $b_1 = 4$ and $b_3 = 8$. Thus, the Betti numbers data do not give us any obstruction for the manifolds to admit a Sasakian metric. However, we prove that, while both nilmanifolds are endowed with a left-invariant $K$-contact structure, they cannot admit any Sasakian metric compatible with the underlying contact structure since they do not satisfy the Hard Lefschetz Theorem.

Concerning the question of finding non-trivial examples of $K$-contact manifolds which cannot carry any Sasakian metric, we mention the very recent paper \cite{hajduktralle2013}, where compact and simply connected examples in dimensions $\geq 11$ are found by using some geometric techniques based on the notion of contact fatness developed by Lerman in \cite{lerman1988} and \cite{lerman2004}. However, we notice that in \cite{hajduktralle2013} the fact that the manifold is non-Sasakian was established by proving that the third Betti number is odd. Finally, when the present paper was in the final stages of preparation, a new result appeared in \cite{Lin2013} where examples of simply connected $K$-contact manifolds without any Sasakian structures in dimensions $\geq 9$ are found by using the Hard Lefschetz Theorem.

\section{Preliminaries}\label{preliminaries}

In this section we recall some basic definitions and properties in contact Riemannian geometry. For further details we refer the reader to the monographs
\cite{blair2010} or \cite{boyergalicki2008}.

Let $M$ be a smooth manifold of dimension $2n+1$. A $1$-form $\eta$ on $M$ is called a \emph{contact form} if $\eta\wedge \d\eta^n$ is a volume form. Then the pair $(M,\eta)$ is called a (strict) \emph{contact manifold}. In any contact manifold one proves the existence of a unique vector field $\xi$, called the \emph{Reeb vector field}, satisfying the properties
\begin{equation*}
i_{\xi}\eta=1, \ \ \ i_{\xi}\d\eta=0.
\end{equation*}
Given a contact manifold $(M,\eta)$ there always exists a Riemannian metric $g$ and a tensor field $\phi$ of type $(1,1)$ such that the following
conditions hold
\begin{gather}
\eta = g(\cdot, \xi),   \label{one}   \\
\d\eta=2g(\cdot,\phi\, \cdot),   \label{two}     \\
\phi^2 = -I + \eta\otimes\xi,   \label{three}
\end{gather}
where $I: TM \to TM$ denotes the identity mapping. From \eqref{one}--\eqref{three} it follows that $\phi\xi
= 0$, $\eta\circ\phi = 0$ and
\begin{equation*}
g(\phi X, \phi Y) = g(X,Y) - \eta(X)\eta(Y)
\end{equation*}
for any $X,Y\in\Gamma(TM)$. Moreover, \eqref{two} implies that the bilinear form $\Phi:=g(\cdot,\phi\cdot)$ is in fact a $2$-form, which is sometimes
called \emph{Sasaki form}.
The manifold $M$ together with the above geometric structure $(\phi,\xi,\eta,g)$ is called a \emph{contact metric manifold}.

A \emph{Sasakian manifold} is a contact metric manifold for which the following normality condition is satisfied
        \begin{equation*}
             \left[ \phi, \phi \right]_{FN} +  2 \d\eta \otimes
            \xi = 0,
        \end{equation*}
            where $[-,-]_{FN}$ is the Fr\"olicher-Nijenhuis
        bracket (see e.g. \cite{kolarmichor}).

An equivalent way to express the Sasakian condition is to say that the Riemannian cone of $(M, \eta,g)$ is a K\"ahler manifold.

In any Sasakian manifold the Reeb vector field is Killing. This last property is equivalent to the condition ${\mathcal L}_{\xi}\phi=0$. More generally, a
contact metric manifold whose Reeb vector field is Killing is called \emph{$K$-contact}. Thus any Sasakian manifold is $K$-contact and it is known that the
converse holds in dimension $3$.

\medskip

A well-known obstruction to the existence of a Sasakian structure on a contact manifold $(M,\eta)$ is given by the following theorem, due to Fujitani.

\begin{theorem}[\cite{fujitani}]\label{bettisasakian}
Let $(M,\eta)$ be a compact contact manifold of dimension $2n+1$. If $M$ admits a compatible Sasakian structure, then for any odd integer $p\leq n+1$ the
Betti numbers $b_p$ are even.
\end{theorem}

In \cite{cappellettidenicolayudin2013} an obstruction stronger than the one expressed in Theorem \ref{bettisasakian} was found, by  proving an odd
dimensional counterpart of the celebrated Hard Lefschetz Theorem.

\begin{theorem}[\cite{cappellettidenicolayudin2013}]\label{HLT}
Let $(M,\phi,\xi,\eta,g)$ be a compact Sasakian manifold of dimension $2n+1$. Then for each integer $0 \leq p \leq n$ the maps
\begin{equation}\label{lef}
 \begin{aligned}
        \lef_{n-p}\colon    H^{n-p}(M) &\to H^{n+p+1}(M)\\
    \left[\, \beta \,\right] &\mapsto \left[\, \eta\wedge (\d\eta)^{p}\wedge\pr\,  \beta \,\right],
\end{aligned}
\end{equation}
are isomorphisms,   $\Pi_\vartriangle  \beta$ denoting the orthogonal projection of $\beta$ on the space of harmonic forms.
\end{theorem}

Notice that, contrary to the even dimensional case, it is not true that the wedge multiplication by $\eta$ or by $\eta\wedge \d\eta$ maps harmonic forms into harmonic forms, so that one is forced to use the metric in order to define the Lefschetz maps \eqref{lef}. Thus, \emph{a priori}, one could expect that
different Sasakian metrics, all compatible with the same underlying contact form,  could lead to different Lefschetz isomorphisms. However, in
\cite[Theorem 4.5]{cappellettidenicolayudin2013}  it is proved that the Lefschetz isomorphisms are independent of the metric. The proof of \cite[Theorem
4.5]{cappellettidenicolayudin2013} suggests to introduce the notion of \emph{contact Lefschetz manifold}. Namely, according to
\cite{cappellettidenicolayudin2013}, let us define the \emph{Lefschetz relation} between cohomology groups $H^p(M)$ and $H^{2n+1-p}(M)$ of a contact
manifold $(M, \eta)$ to be
\begin{equation*}
    \rel_{\lef_p} = \left\{\, \left(\,\, [ \beta ]\,\,,\, [ \epseta L^{n-p}\beta
    ] \,\,\right) \,\middle|\, \beta \in \Omega^p(M),\ \d\beta = 0,\ \ixi\beta =0,\
    L^{n-p+1}\beta =0\, \right\}.
\end{equation*}
Here, given $\alpha\in \Omega^k(M)$, the operator $\eps_\alpha$ is defined by
$$\eps_\alpha \beta=\alpha \wedge\beta,$$
$L=\eps_{\frac12\d\eta}$ and the conditions $\d\beta = 0$, $\ixi\beta =0$, and $L^{n-p+1}\beta =0$ imply that the $((2n+1)-p)$-form $\epseta L^{n-p}\beta$ is closed and $[\epseta L^{n-p}\beta]=\lef_p[\beta]$ (see Theorem 4.5 in \cite{cappellettidenicolayudin2013}).
Therefore if $(M,\eta)$ admits a compatible Sasakian metric,  due to Theorem~\ref{HLT} it follows that $\rel_{\lef_p}$ is the graph of the isomorphism
$\lef_p$. More generally, a compact contact manifold $(M,\eta)$ is said to be \emph{Lefschetz contact}  if it satisfies the \emph{hard Lefschetz property}, that is, for every $p\le n$ the relation $\rel_{\lef_p}$ is the graph of an isomorphism between $H^p(M)$ and $H^{2n+1-p}(M)$.

We point out that, according to \cite[Theorem 5.2]{cappellettidenicolayudin2013}, every Lefschetz contact manifold satisfies the restrictions on the Betti numbers for compact Sasakian manifolds stated in Theorem \ref{bettisasakian}.

\section{Examples of compact $K$-contact manifolds with no compatible Sasakian metric}

Hereinafter, by using the aforementioned Hard Lefschetz Theorem for Sasakian manifolds, we shall present two examples of compact contact manifolds
(actually, $K$-contact manifolds) of dimension $5$ and $7$, respectively, which do not admit any compatible Sasakian metric. As we shall see, one cannot deduce the assertion from the well-known restrictions concerning the Betti numbers of the manifold.

\medskip

\begin{ex}\label{example5D}
Let us consider the 5-dimensional nilpotent Lie algebra $\mathfrak{g}$ with non-zero Lie brackets
\begin{gather}\label{lie-brackets5}
[X_{1},X_{2}]=X_{3}, \ \ \ [X_{1},X_{3}]=X_{4}, \ \ \ [X_{1},X_{4}]=X_{5}, \ \ \ [X_{2},X_{3}]=X_{5}.
\end{gather}

Let us denote by $\alpha_i$ the dual $1$-form of the vector $X_i$, for any $i\in\left\{1,\ldots,5\right\}$. Then by \eqref{lie-brackets5} it follows that
\begin{align}
\begin{split}
\d \alpha_{1}&=0, \ \ \ \d \alpha_{2}=0, \quad \d \alpha_{3}=-\alpha_{1}\wedge \alpha_{2} , \quad \d \alpha_{4}=-\alpha_{1}\wedge \alpha_{3},\\
\d \alpha_{5} &= \alpha_{1}\wedge \alpha_{4} + \alpha_{2}\wedge \alpha_{3}. \label{differentials5}
\end{split}
\end{align}
Here, $\d\colon \wedge^{k}\g^{\ast}\to \wedge^{k+1}\g^{\ast}$ is the Chevalley-Eilenberg differential in the Lie algebra $\g$. We have $\alpha_{5}\wedge (\d\alpha_{5})^2=2 \alpha_{1}\wedge \alpha_{2} \wedge \alpha_{3}\wedge \alpha_{4}\wedge \alpha_{5}$, hence  if $\eta$ is the  left-invariant 1-form induced by $\alpha_5$ on the simply connected Lie group $G$ whose Lie algebra is $\mathfrak{g}$, then $\eta$ is a contact 1-form on $G$. The corresponding Reeb vector field is the left-invariant vector field $\xi$ on $G$ induced by $X_5$.

Now, let us consider a co-compact discrete subgroup $\Gamma$ of $G$. Such a subgroup exists by the Malcev's criterion, since in our case the structure constants are  integers. Then the quotient $M=G/\Gamma$ is a compact nilmanifold.
\begin{rmk}
The Lie group $G$ is isomorphic to $\Rr^5$ endowed with the multiplication defined by
\begin{equation*}
\begin{split}
(x^1,x^2,x^3,x^4,x^5)\cdot&(y^1,y^2,y^3,y^4,y^5)= \left(x^1+y^1, x^2+y^2, x^3+y^3+x^1 y^2,x^4+y^4\right.\\
 & {}+x^1 y^3+\frac{(x^1)^2}{2}y^2,x^5+y^5-x^1 y^4-\left(\frac{(x^1)^2}{2}+x^2\right)y^3\\
 &  {}- \frac{x^1}{2} (y^2)^2-\left.\left(\frac{(x^1)^3}{6}+x^1x^2\right)y^2\right),
\end{split}
\end{equation*}
for $(x^1,x^2,x^3,x^4,x^5),(y^1,y^2,y^3,y^4,y^5)\in\Rr^5$. A basis of the space of the left-invariant 1-forms is given by
\begin{align*}
\begin{split}
\alpha_{1}&=\d x^1, \quad \alpha_{2}=\d x^2, \quad \alpha_{3}=\d x^3-x^1\d x^2, \quad\alpha_{4}=\d x^4-x^1\d x^3+ \frac{(x^1)^2}{2}\d x^2,\label{basisleft}\\
\alpha_{5}&= \d x^5+x^1\d x^4- \frac{(x^1)^2}{2}\d x^3+\frac{(x^1)^3}{6}\d x^2+x^2\d x^3.
\end{split}
\end{align*}
As a co-compact discrete subgroup of $G$, we can take for example
\begin{equation*}
\Gamma= 6\Z \times \Z^4.
\end{equation*}
\end{rmk}
By the Nomizu's Theorem, $H^{\ast}_{DR}(M)\cong H^{\ast}(\mathfrak{g})$,  where
$H^{\ast}(\mathfrak{g}) = H^{\ast}(\mathfrak{g},\R)$ denotes the Chevalley-Eilenberg cohomology of the Lie algebra
$\mathfrak{g}$ with the coefficients in the trivial module $\R$.

By \eqref{differentials5} we have that
\begin{gather*}
H^{1}(\mathfrak{g})=\langle [\alpha_{1}],  [\alpha_{2}] \rangle.
\end{gather*}
Thus $b_{1}(M)=2$, so that the Betti numbers do not give any obstruction for $(M,\eta)$ to admit a Sasakian structure.  Now, according to \cite[Section 5]{cappellettidenicolayudin2013},
let us consider the Lefschetz relation
\begin{equation*}
{\mathcal R}_{Lef_p}=\left\{([\beta],[\eta\wedge(\d\eta)^{2-p}\wedge\beta]) \  | \ \beta\in\Omega^{p}(M), \ \d\beta=0, \ i_{\xi}\beta=0, \
(\d\eta)^{3-p}\wedge\beta=0\right\}.
\end{equation*}
As a consequence of \cite[Theorem 4.5]{cappellettidenicolayudin2013}, if a $5$-dimensional contact manifold $(M,\eta)$ admits any compatible Sasakian metric then ${\mathcal R}_{Lef_p}$ is the graph of an
isomorphism for every $p\leq 2$. Actually, let us prove that in our case the above property does not hold for $p=1$. Let us consider the  $1$-form
$\beta$ on $M$ induced by the left-invariant 1-form on $G$ associated to $\alpha_{2}$. Clearly $\beta$ is closed, $i_\xi \beta=0$ and
\begin{equation*}
(\d\eta)^2\wedge\beta=0.
\end{equation*}
Moreover, using that $[\alpha_{2}]$ is a non-zero element of $H^{1}(\mathfrak{g})$, we deduce that $[\beta]$ is a non-zero element of $H^{1}_{DR}(M)$.  We prove that
\begin{equation*}
[\eta\wedge \d\eta\wedge\beta]=0.
\end{equation*}
Indeed we have
\begin{equation*}
\d (\alpha_3 \wedge \alpha_4 \wedge \alpha_5) =-\alpha_{1}\wedge \alpha_{2}\wedge \alpha_{4}\wedge \alpha_{5}.
\end{equation*}
Therefore, $\eta \wedge \d \eta \wedge \beta$ is an exact 4-form on $M$. In fact, $\eta \wedge \d \eta \wedge \beta=-\d \gamma$, with $\gamma$ the 3-form on $M$ induced by the left-invariant 3-form on $G$ associated to $\alpha_3 \wedge \alpha_4 \wedge \alpha_5\in\wedge^{3}\g^{\ast}$. So, ${\mathcal R}_{Lef_1}$ is not the graph of an isomorphism and we conclude that $M$ cannot carry any (not necessarily left-invariant) Sasakian
metric compatible with $\eta$.

\begin{rmk}
We point out that though $(M,\eta)$ cannot admit a Sasakian metric, it admits a (even left-invariant) $K$-contact structure. Indeed, we can define an
endomorphism $\phi:\mathfrak{g}\longrightarrow \mathfrak{g}$ by setting
\begin{equation*}
\phi X_1 =-X_4, \ \ \ \phi X_2 =- X_3, \ \ \ \phi X_3 = X_2, \ \ \ \phi X_4 = X_1, \ \ \ \phi X_5 =0
\end{equation*}
and we can define a positive definite bilinear form $g$ by declaring that $\left\{X_{1},\ldots,X_{5}\right\}$ is $g$-orthonormal. Then
$(\phi,\xi,\eta,g)$ induces a left-invariant contact metric structure on $G$ which descends to the quotient. Since $X_{5}$ belongs to the center of the Lie algebra $\mathfrak{g}$, we have immediately that the Reeb vector field $\xi$ is Killing and thus the structure is $K$-contact.
\end{rmk}
\end{ex}

\medskip

\begin{ex}\label{example7D}
Now, we shall present an example of a compact $7$-dimensional $K$-contact manifold which does not admit any compatible Sasakian metric. As we shall see, one cannot deduce this fact by Theorem \ref{bettisasakian}.

\medskip

In a recent paper \cite{kutsak2012}, Kutsak has classified all nilpotent $7$-dimensional Lie algebras carrying a left-invariant contact structure. Let us consider the example (1457B) in that paper, namely the nilpotent Lie algebra $\mathfrak{g}$ with non-zero Lie brackets
\begin{align}
\begin{split}\label{lie-brackets}
[X_{1},X_{2}]&=X_{3}, \qquad [X_{1},X_{3}]=X_{4}, \qquad  [X_{1},X_{4}]=X_{7}, \\
[X_{2},X_{3}]&=X_{7}, \qquad  [X_{5},X_{6}]=X_{7}.
\end{split}
\end{align}
Let us denote by $\alpha_i$ the dual $1$-form of the vector $X_i$, for any $i\in\left\{1,\ldots,7\right\}$. Then by \eqref{lie-brackets} it follows that
\begin{align}
\begin{split}
\d \alpha_{1}&=0, \qquad \d \alpha_{2}=0, \qquad \d \alpha_{5}=0, \qquad \d \alpha_{6}=0, \qquad \d \alpha_{3}=-\alpha_{1}\wedge \alpha_{2}, \label{differentialss}\\
\d\alpha_{4}&=-\alpha_{1}\wedge \alpha_{3} \qquad \d \alpha_{7} = -\alpha_{1}\wedge \alpha_{4} - \alpha_{2}\wedge \alpha_{3} - \alpha_{5}\wedge \alpha_{6}.
\end{split}
\end{align}
Let $G$ be the simply connected Lie group $G$ with Lie algebra $\mathfrak{g}$ and consider the left-invariant 1-form on $G$ associated to $\alpha_{7}$. Then, such a 1-form defines a contact structure on $G$. The corresponding Reeb vector field $\xi$ is the left-invariant vector field on $G$ induced by $X_7$.

Now, let us consider a co-compact discrete subgroup $\Gamma$ of $G$. Then the quotient $M=G/\Gamma$ is a compact nilmanifold.
\begin{rmk}
The Lie group $G$ is isomorphic to $\Rr^7$ endowed with the multiplication defined by
\begin{equation*}
\begin{split}
(x^1,x^2,x^3,x^4,x^5,x^6,x^7)&\cdot(y^1,y^2,y^3,y^4,y^5,y^6,y^7)= \left(x^1+y^1, x^2+y^2, x^3+y^3+x^1 y^2,\right.\\
 & x^4+y^4+x^1 y^3+\frac{(x^1)^2}{2}y^2,x^5+y^5,x^6+y^6, x^7+y^7+x^1 y^4\\
 &  {}+ \left(\frac{(x^1)^2}{2}+x^2\right)y^3+\frac{x^1}{2} (y^2)^2 +\left.\left(\frac{(x^1)^3}{6}+x^1x^2\right)y^2 +x^5y^6\right),
\end{split}
\end{equation*}
for $(x^1,x^2,x^3,x^4,x^5,x^6,x^7),(y^1,y^2,y^3,y^4,y^5,y^6,y^7)\in\Rr^7$. A basis of the space of the left-invariant 1-forms is given by
\begin{align*}
\begin{split}
\alpha_{1}&=\d x^1, \quad \alpha_{2}=\d x^2, \quad \alpha_{3}=\d x^3-x^1\d x^2, \quad\alpha_{4}=\d x^4-x^1\d x^3+ \frac{(x^1)^2}{2}\d x^2,\label{basisleft}\\
\alpha_{5}&= \d x^5, \,\alpha_{6}= \d x^6,\, \alpha_{7}= \d x^7 -x^1\d x^4+ \frac{(x^1)^2}{2}\d x^3-\frac{(x^1)^3}{6}\d x^2-x^2\d x^3-x^5\d x^6.
\end{split}
\end{align*}
A co-compact discrete subgroup of $G$ is
\begin{equation*}
\Gamma= 6\Z \times \Z^6.
\end{equation*}
\end{rmk}

By the Nomizu's Theorem $H^{\ast}_{DR}(M)\cong H^{\ast}(\mathfrak{g})$. By \eqref{differentialss} we have that
$
H^{1}(\mathfrak{g})=\langle [\alpha_{1}],  [\alpha_{2}], [\alpha_{5}],
[\alpha_{6}] \rangle$.
Thus $b_{1}(M)=4$. Furthermore, by using the library~\cite{plural} of the
computer algebra system \texttt{Singular}~\cite{singular}, one can compute the
dimension of $H^{3}(\mathfrak{g})$. Namely, the cohomology group
$H^3(\mathfrak{g},\R)$ is isomorphic to the third
Ext-group $\Ext_{\mathfrak{U}(\mathfrak{g})}^3(\R,\R)$, where
$\mathfrak{U}(\mathfrak{g})$ denotes  the universal enveloping algebra
of $\mathfrak{g}$.
Note that the Lie algebra $\mathfrak{g}$ is a positively graded Lie algebra with the grading
\begin{align*}
\deg(X_1)&= \deg(X_5)=1, & \deg(X_2) &= 2, & \deg(X_3)&=3, \\ \deg(X_4)&=
\deg(X_6) = 4, & \deg(X_7) &= 5.
\end{align*}
Therefore, we can consider $\mathfrak{U}(\mathfrak{g})$ as a connected graded
associative algebra.
Thus the dimension of
$\Ext_{\mathfrak{U}(\mathfrak{g})}^k(\R,\R)$ coincides with the rank of the
(necessarily free) $k$-th module of the minimal projective resolution of $\R$ over
$\mathfrak{U}(\mathfrak{g})$. One can use the procedure \texttt{mres} of
\texttt{Singular} to construct such a resolution. As a result, we get that
$b_{3}(M)=8$. Therefore, the dimension of $H^{3}(\mathfrak{g})$ does not give any obstruction to the existence of a Sasakian structure on $(M,\eta)$, where $\eta$  is the contact 1-form on $M$ induced by the left-invariant 1-form on $G$ associated to $\alpha_7$.

Now, according to \cite[Section 5]{cappellettidenicolayudin2013},
let us consider the Lefschetz relation
\begin{equation*}
{\mathcal R}_{Lef_p}=\left\{([\beta],[\eta\wedge(\d\eta)^{3-p}\wedge\beta]) \  | \ \beta\in\Omega^{p}(M), \ \d\beta=0, \ i_{\xi}\beta=0, \
(\d\eta)^{4-p}\wedge\beta=0\right\}.
\end{equation*}
As recalled in Section \ref{preliminaries}, as a consequence of \cite[Theorem 4.5]{cappellettidenicolayudin2013}, if a $7$-dimensional contact manifold
$(M,\eta)$ admits any compatible Sasakian metric then ${\mathcal R}_{Lef_p}$ is the graph of an isomorphism for every $p\leq 3$. Actually, let us prove
that in our case the above property does not hold for $p=1$. Let us consider the  $1$-form $\beta$ on $M$ induced by the left-invariant 1-form on $G$
associated to $\alpha_{1}\in \mathfrak{g}^\ast$. First of all let us check that ${\mathcal R}_{Lef_1}$ is defined for $\beta$. Indeed, $\beta$ is
closed, $i_\xi \beta=0$ and
\begin{equation*}
(\d\eta)^{3}\wedge\beta=0.
\end{equation*}
Clearly, $[\beta]$ defines a non-zero element of $H^{1}_{DR}(M)$. Finally
\begin{equation*}
[\eta\wedge (\d\eta)^2\wedge\beta]=0.
\end{equation*}
Indeed we have
\begin{align*}
\alpha_{7}\wedge (\d \alpha_{7})^2 \wedge \alpha_{1}&=- 2 \alpha_{1}\wedge \alpha_{2}\wedge \alpha_{3}\wedge \alpha_{5}\wedge\alpha_{6}\wedge \alpha_{7}\\
&=2\d (\alpha_{2} \wedge\alpha_{4} \wedge \alpha_{5}\wedge\alpha_{6}\wedge \alpha_{7}).
\end{align*}
Therefore ${\mathcal R}_{Lef_1}$ is not the graph of an isomorphism and we conclude that $M$ cannot carry any (not necessarily left-invariant) Sasakian
metric compatible with $\eta$.

\begin{rmk}
We point out that though $(M,\eta)$ cannot admit a Sasakian metric, it admits a (even left-invariant) $K$-contact structure. Indeed, we can define an
endomorphism $\phi:\mathfrak{g}\longrightarrow \mathfrak{g}$ by setting
\begin{equation*}
\begin{split}
\phi X_1 &= X_4, \qquad \phi X_2 = X_3,  \qquad \phi X_3 = -X_2, \qquad \phi X_4 = -X_1,\\
\phi X_5 &= X_6, \qquad \phi X_6 = -X_5, \qquad \phi X_7 =0.
\end{split}
\end{equation*}
and we can define a positive definite bilinear form $g$ by declaring that $\left\{X_{1},\ldots,X_{7}\right\}$ is $g$-orthonormal. Then
$(\phi,X_7,\alpha_{7},g)$ induces a left-invariant contact metric structure on $G$ which descends to the quotient. Since $X_7$ belongs to the
center of the Lie algebra $\mathfrak{g}$, we have immediately that the Reeb vector field is Killing and thus the structure is $K$-contact.
\end{rmk}
\end{ex}

\begin{rmk}
In \cite{andradafinovezzoni2009} it was proved that every nilpotent Lie algebra admitting a Sasakian structure is necessarily
isomorphic to the Heisenberg Lie algebra. Thus one could think that the non-Sasakian property for our examples can be  deduced from this result. Actually,
our results are stronger, because we prove that the nilmanifolds in Examples \ref{example5D} and  \ref{example7D} cannot admit \emph{any} compatible
Sasakian structure, not necessarily left-invariant.
\end{rmk}

\medskip

\noindent{\bf Acknowledgments}

\noindent Research partially supported by CMUC, funded by the European program
COMPETE/FEDER, by FCT (Portugal) grants PEst-C/MAT/UI0324/2011 (A.D.N. and I.Y.), by MICINN (Spain) grants
MTM2011-15725-E, MTM2012-34478 (A.D.N. and J.C.M.), the project of the Canary Government ProdID20100210 (J.C.M.), and by Prin 2010/11 -- Variet\`{a} reali e complesse: geometria, topologia e analisi armonica –- Italy (B.C.M.).


\begin{thebibliography}{99}
\bibitem{andradafinovezzoni2009} A. Andrada, A. Fino, L. Vezzoni, \emph{A class of Sasakian 5-manifolds}, Transformation Groups \textbf{14} (2009),
493--512.

\bibitem{bensongordon} C. Benson, C.~S. Gordon, \emph{K\"{a}hler and symplectic structures on nilmanifolds}, Topology \textbf{27} (1988), 513--518.

\bibitem{blair2010} D.~E. Blair, \emph{Riemannian geometry of contact and symplectic manifolds}, Second Edition. Progress in Mathematics \textbf{203}, Birkh\"{a}user, Boston, 2010.

\bibitem{boyergalicki2008} C.~P. Boyer and K. Galicki, \emph{Sasakian geometry}, Oxford University Press, 2008.

\bibitem{cappellettidenicolayudin2013}B. Cappelletti-Montano, A. De Nicola, I. Yudin, \emph{Hard Lefschetz Theorem for Sasakian manifolds}, arXiv:1306.2896v1.

\bibitem{chineadeleonmarrero}
D.~Chinea, M.~de~Le{\'o}n, and J.~C. Marrero, \emph{Topology of cosymplectic
  manifolds}, J. Math. Pures Appl. (9) \textbf{72} (1993), no.~6, 567--591.

\bibitem{singular}
	W. Decker, G.-M.~Greuel, G.~Pfister, H.~Sch{\"o}nemann,
	\newblock {\sc Singular} {3-1-5} --- \emph{{A} computer algebra system for polynomial computations},
	\newblock {http://www.singular.uni-kl.de}, 2012.


\bibitem{fujitani} T. Fujitani, \emph{Complex-valued differential forms on normal contact Riemannian manifolds},
T\^{o}hoku Math. J. \textbf{18} (1966), 349–-361.

\bibitem{hajduktralle2013} B. Hajduk, A. Tralle, \emph{On simply connected $K$-contact non-Sasakian manifolds}, arXiv:1305.2591v1.

\bibitem{kolarmichor}I. Kol{\'a}{\v{r}}, P. Michor, J. Slov{\'a}k, \emph{Natural operations in differential geometry}, Springer-Verlag, Berlin, 1993.

\bibitem{kutsak2012} S. Kutsak, \emph{Invariant contact structures on 7-dimensional nilmanifolds}, Geom. Dedicata \textbf{172} (2014), 351--361.

\bibitem{lerman1988} E. Lerman, \emph{How fat is a fat bundle?}, Lett. Math. Phys. \textbf{15} (1988), 335--339.

\bibitem{lerman2004} E. Lerman, \emph{Contact fiber bundles}, J. Geom. Phys. \textbf{49} (2004), 52--66.

\bibitem{plural} V. Levandovskyy, F.~L.~Lobillo, C. Rabelo, O. Motsak,
{\tt nctools.lib}. {A} {\sc Singular} {3-1-5}  library of general tools for noncommutative algebras, 2013.

\bibitem{Lin2013} Y. Lin, \emph{Lefschetz contact manifolds and odd dimensional symplectic geometry}, arXiv:1311.1431v2.
	
\bibitem{thurston} W.~P. Thurston, \emph{Some simple examples of symplectic manifolds}, Proc. Amer. Math. Soc., \textbf{55} (1976), 467--468.

\bibitem{tralleoprea1997}A. Tralle, G. Oprea, \emph{Symplectic manifolds with no K\"{a}hler structure}, Lecture Notes in Mathematics 1661, Springer-Verlag, Berlin, 1997.

\end{thebibliography}
\end{document}